New series representations of $\pi$, $\pi^3$ and $\pi^5$ in terms of Euler numbers and $\pi^2$, $\pi^4$ and $\pi^6$ in terms of Bernoulli numbers

H.C.Gupta

Physics Department

Indian Institute of Technology

New Delhi-110016, INDIA

hcgupta@physics.iitd.ac.in

ABSTRACT

New series representations for odd powers of $\pi$ i.e. $\pi$, $\pi^3$ and $\pi^5$ in terms of Euler numbers and even powers of $\pi$ i.e. $\pi^2$, $\pi^4$ and $\pi^6$ in terms of Bernoulli numbers have been obtained empirically.

Keywords: Powers of $\pi$, Series representation, Bernoulli numbers, Euler numbers

Mathematics Subject Classification 2010: 11B68, 11B50

## 1. Introduction

The representation of the mathematical constants π and its powers by series, products and integrals have attracted the attention of researchers for a long time. The classic result for $\pi$, $\pi^2$, $\pi^3$, $\pi^4$, $\pi^5$, $\pi^6$ in terms of a series of powers of integers was presented by Leibnitz and Euler as

$$\pi = 4 \sum_{k=0}^{\infty} \frac{(-1)^k}{(2k+1)} \tag{1}$$

$$\pi^2 = 6 \sum_{k=1}^{\infty} \frac{1}{k^2} \tag{2}$$

$$\pi^3 = 32 \sum_{k=0}^{\infty} \frac{(-1)^k}{(2k+1)^3} \quad , \quad \pi^4 = 90 \sum_{k=1}^{\infty} \frac{1}{k^4} \tag{3}$$

$$\pi^5 = \frac{1536}{5} \sum_{k=0}^{\infty} \frac{(-1)^k}{(2k+1)^5} \quad , \quad \pi^6 = 945 \sum_{k=1}^{\infty} \frac{1}{k^6} \tag{4}$$

Borwein[5] had described very elegantly the history of π along with its selected useful formulas in terms of series representation. Alzer and Koumandos[2] presented a series representation for π as

$$\pi = 4 \sum_{k=0}^{\infty} \frac{1}{(1+\mu)^{k+1}} \sum_{m=0}^{k} \binom{k}{m} \frac{(-1)^m \mu^{k-m}}{(2m+1)} \tag{5}$$

where they introduced a positive parameter μ, on the lines of modifying the Amore's technique[3] to obtain one-parameter family of series representation of mathematical constants. Further, in case of mathematical constants involving π, series representations were given by Alzer etal.[1] in terms of digamma function. Recently, Blagouchine[4] obtained new series expansions for generalized Euler constants containing polynomials in $\pi^{-2}$ with rational coefficients.

Alzer etal.[1] gave series representation for $\pi^2$ as

$$\pi^2 = 4 \sum_{k=1}^{\infty} \frac{\mu_k h_k}{k} \tag{6}$$

where $\mu_n$ denotes the normalized binomial mid-coefficient defined by

$$\mu_k = \frac{1.3.5.........(2k-1)}{2.4.6. .........2k} \tag{7}$$

and
$$h_k = \sum_{n=1}^{k} \frac{1}{(2n-1)} \tag{8}$$

Kolbig [6] gave an interesting relation as

$$\pi^2 = 2 \sum_{k=1}^{\infty} \frac{\sigma_k}{k} \tag{9}$$

where
$$\sigma_n = p_n \sum_{k=1}^{n} \frac{1}{(4k-1)} + q_n \sum_{k=1}^{n} \frac{1}{(4k-3)} \tag{10}$$

$$p_n = \prod_{k=1}^{n} \frac{(4k-1)}{4k} \quad, \quad q_n = \prod_{k=1}^{n} \frac{(4k-3)}{4k} \tag{11}$$

Yet another series representation for $\pi^2$ were given by Alzer etal.[1] as

$$\pi^2 = 3 \sum_{k=1}^{\infty} \frac{\mu_k H_k}{k} \tag{12}$$

where
$$H_n = \sum_{k=1}^{n} \frac{1}{k} \tag{13}$$

In this paper, the series representations for $\pi$, $\pi^3$ and $\pi^5$ are presented where the series involves the application of Euler numbers $E_{2k}$, defined as

$$\sum_{m=1}^{\infty} \frac{(-1)^{m+1}}{(2m-1)^{2k+1}} = \frac{\pi^{2k+1}}{2^{2k+2} (2k)!} |E_{2k}| \tag{14}$$

with $E_0 = 1$ giving $E_2 = -1$, $E_4 = 5$, $E_6 = -61$, $E_8 = 1385$, $E_{10} = -50521$, $E_{12} = 2702765$, …… and the series $\sum_{j=0}^{k-1} \left[ -\frac{1}{(2m-1)^2 \pi^2} \right]^j \frac{1}{(2k-2j-1)!}$ for integers j, k ≥ 0.

Further, the series representations for $\pi^2$, $\pi^4$ and $\pi^6$ are presented where the series involves the application of Bernoulli numbers $B_{2k}$, defined as

$$\sum_{m=1}^{\infty} \frac{1}{(m)^{2k}} = (-1)^{k-1} \frac{(2\pi)^{2k}}{2 (2k)!} B_{2k}, \ k = 1, 2, \ldots \tag{15}$$

with $B_0 = 1$ giving $B_1 = -\frac{1}{2}$, $B_2 = \frac{1}{6}$, $B_4 = -\frac{1}{30}$, $B_6 = \frac{1}{42}$, $B_8 = -\frac{1}{30}$, $B_{10} = \frac{5}{66}$, $B_{12} = -\frac{691}{2730}$, …… and the series $\sum_{j=0}^{k-1} \left[ -\frac{1}{(m)^2 \pi^2} \right]^j \frac{1}{(2k-2j-1)!}$ for integers j, k ≥ 0.

## 2. Series representations for odd powers $\pi$, $\pi^3$ and $\pi^5$

### 2a. The series for $\pi$

The series for $\pi$ is presented as

$$\pi = \sum_{n=1}^{\infty} \frac{(-1)^{n+1}}{(2n-1)} 2^{2k+2} \ (2k+1)! \ \sum_{j=0}^{k} \left[-\frac{1}{(2n-1)^2 \pi^2}\right]^j \frac{1}{(2k-2j+1)!} \qquad (16)$$

where j, k are integers and $\geq 0$.

**Verification of the series for some k values**

$k = 0$ gives series as $\sum_{n=1}^{\infty} \frac{(-1)^{n+1}}{(2n-1)} 2^2 \ (1)! \ \left[-\frac{1}{(2n-1)^2 \pi^2}\right]^0 \frac{1}{(1)!} = 4 \cdot \frac{\pi}{4} = \pi$

$k = 1$ gives series as $\sum_{n=1}^{\infty} \frac{(-1)^{n+1}}{(2n-1)} 2^4 \ (3)! \ \left[\frac{1}{3!} - \frac{1}{(2n-1)^2 \pi^2}\right] = 16 \left[\frac{\pi}{4} - \frac{6}{\pi^2} \frac{\pi^3}{32}\right] = \pi$

$k = 2$ gives series as $\sum_{n=1}^{\infty} \frac{(-1)^{n+1}}{(2n-1)} 2^6 \ (5)! \ \left[\frac{1}{5!} - \frac{1}{3!} \frac{1}{(2n-1)^2 \pi^2} + \frac{1}{(2n-1)^4 \pi^4}\right] = 64 \left[\frac{\pi}{4} - \frac{20}{\pi^2} \frac{\pi^3}{32} + \frac{120}{\pi^4} \frac{5\pi^5}{1536}\right] = \pi$

$k = 3$ gives series as $\sum_{n=1}^{\infty} \frac{(-1)^{n+1}}{(2n-1)} 2^8 \ (7)! \ \left[\frac{1}{7!} - \frac{1}{5!} \frac{1}{(2n-1)^2 \pi^2} + \frac{1}{3!} \frac{1}{(2n-1)^4 \pi^4} - \frac{1}{(2n-1)^6 \pi^6}\right] =$

$256 \left[\frac{\pi}{4} - \frac{42}{\pi^2} \frac{\pi^3}{32} + \frac{840}{\pi^4} \frac{5\pi^5}{1536} - \frac{5040}{\pi^6} \frac{61\pi^7}{184320}\right] = \pi$

$k = 4$ gives series as $\sum_{n=1}^{\infty} \frac{(-1)^{n+1}}{(2n-1)} 2^{10} \ (9)! \ \left[\frac{1}{9!} - \frac{1}{7!} \frac{1}{(2n-1)^2 \pi^2} + \frac{1}{5!} \frac{1}{(2n-1)^4 \pi^4} - \frac{1}{3!} \frac{1}{(2n-1)^6 \pi^6} + \frac{1}{(2n-1)^8 \pi^8}\right] =$

$1024 \left[\frac{\pi}{4} - \frac{72}{\pi^2} \frac{\pi^3}{32} + \frac{3024}{\pi^4} \frac{5\pi^5}{1536} - \frac{60480}{\pi^6} \frac{61\pi^7}{184320} + \frac{362880}{\pi^8} \frac{1385\pi^9}{41287680}\right] = \pi$

On similar lines, the series representation (16) can be proved for any integer k value.

### 2b. The series for $\pi^3$

The series for $\pi^3$ is presented as

$$\pi^3 = \sum_{n=1}^{\infty} \frac{(-1)^{n+1}}{(2n-1)^3 (2^{2k+2}-1)} 2^{2k+4} \ (2k+3)! \ \sum_{j=0}^{k} \left[-\frac{1}{(2n-1)^2 \pi^2}\right]^j \frac{1}{(2k-2j+1)!} \qquad (17)$$

where j, k are integers and $\geq 0$.

**Verification of the series for some k values**

k = 0 gives series as $\sum_{n=1}^{\infty} \frac{(-1)^{n+1}}{(2n-1)^3} 2^4 \frac{1}{(2^2-1)}$ (3)! $\left[-\frac{1}{(2n-1)^2 \pi^2}\right]^0 \frac{1}{(1)!} = 32 \cdot \frac{\pi^3}{32} = \pi^3$

k = 1 gives series as $\sum_{n=1}^{\infty} \frac{(-1)^{n+1}}{(2n-1)^3} 2^6 \frac{1}{(2^4-1)}$ (5)! $\left[\frac{1}{3!} - \frac{1}{(2n-1)^2 \pi^2}\right] = 512 \left[\frac{1}{6} \frac{\pi^3}{32} - \frac{1}{\pi^2} \frac{5\pi^5}{1536}\right] = \pi^3$

k = 2 gives series as $\sum_{n=1}^{\infty} \frac{(-1)^{n+1}}{(2n-1)^3} 2^8 \frac{1}{(2^6-1)}$ (7)! $\left[\frac{1}{5!} - \frac{1}{3!} \frac{1}{(2n-1)^2 \pi^2} + \frac{1}{(2n-1)^4 \pi^4}\right] =$

$20480 \left[\frac{1}{120} \frac{\pi^3}{32} - \frac{1}{6\pi^2} \frac{5\pi^5}{1536} + \frac{1}{\pi^4} \frac{61\pi^7}{184320}\right] = \pi^3$

k = 3 gives series as $\sum_{n=1}^{\infty} \frac{(-1)^{n+1}}{(2n-1)^3} 2^{10} \frac{1}{(2^8-1)}$ (9)! $\left[\frac{1}{7!} - \frac{1}{5!} \frac{1}{(2n-1)^2 \pi^2} + \frac{1}{3!} \frac{1}{(2n-1)^4 \pi^4} - \frac{1}{(2n-1)^6 \pi^6}\right] =$

$\frac{371589120}{255} \left[\frac{1}{5040} \frac{\pi^3}{32} - \frac{1}{120\pi^2} \frac{5\pi^5}{1536} + \frac{1}{6\pi^4} \frac{61\pi^7}{184320} - \frac{1}{\pi^6} \frac{1385\pi^9}{41287680}\right] = \pi^3$

k = 4 gives series as $\sum_{n=1}^{\infty} \frac{(-1)^{n+1}}{(2n-1)^3} 2^{12} \frac{1}{(2^{10}-1)}$ (11)! $\left[\frac{1}{9!} - \frac{1}{7!} \frac{1}{(2n-1)^2 \pi^2} + \frac{1}{5!} \frac{1}{(2n-1)^4 \pi^4} - \frac{1}{3!} \frac{1}{(2n-1)^6 \pi^6} + \frac{1}{(2n-1)^8 \pi^8}\right] = \frac{4954521600}{31} \left[\frac{1}{362880} \frac{\pi^3}{32} - \frac{1}{5040\pi^2} \frac{5\pi^5}{1536} + \frac{1}{120\pi^4} \frac{61\pi^7}{184320} - \frac{1}{6\pi^6} \frac{1385\pi^9}{41287680} + \frac{1}{\pi^8} \frac{50521\pi^{11}}{14863564800}\right] = \pi^3$

On similar lines, the series representation (17) can be proved for any integer k value.

### 2c. The series for $\pi^5$

The series for $\pi^5$ is presented as

$$\pi^5 = \sum_{n=1}^{\infty} \frac{(-1)^{n+1}}{(2n-1)^5 \{[2^{2k+2}(2k^2+9k+6)]+1\}} 2^{2k+6} \ (2k+5)! \sum_{j=0}^{k} \left[-\frac{1}{(2n-1)^2 \pi^2}\right]^j \frac{1}{(2k-2j+1)!} \qquad (18)$$

where j, k are integers and $\geq 0$.

**Verification of the series for some k values**

k = 0 gives series as $\sum_{n=1}^{\infty} \frac{(-1)^{n+1}}{(2n-1)^5} 2^6 \frac{1}{\{[4(6)]+1\}}$ (5)! $\left[-\frac{1}{(2n-1)^2 \pi^2}\right]^0 \frac{1}{(1)!} = \frac{1536}{5} \cdot \frac{5\pi^5}{1536} = \pi^5$

k = 1 gives series as $\sum_{n=1}^{\infty} \frac{(-1)^{n+1}}{(2n-1)^5} 2^8 \frac{1}{\{[16(17)]+1\}}$ (7)! $\left[\frac{1}{3!} - \frac{1}{(2n-1)^2 \pi^2}\right] = \frac{256(5040)}{273} \left[\frac{1}{6} \frac{5\pi^5}{1536} - \frac{1}{\pi^2} \frac{61\pi^7}{184320}\right] = \pi^5$

k = 2 gives series as $\sum_{n=1}^{\infty} \frac{(-1)^{n+1}}{(2n-1)^5} 2^{10} \frac{1}{\{[64(32)]+1\}}$ (9)! $\left[\frac{1}{5!} - \frac{1}{3!} \frac{1}{(2n-1)^2 \pi^2} + \frac{1}{(2n-1)^4 \pi^4}\right] =$

$\frac{1024(362880)}{2049} \left[\frac{1}{120} \frac{5\pi^5}{1536} - \frac{1}{6\pi^2} \frac{61\pi^7}{184320} + \frac{1}{\pi^4} \frac{1385\pi^9}{41287680}\right] = \pi^5$

k = 3 gives series as $\sum_{n=1}^{\infty} \frac{(-1)^{n+1}}{(2n-1)^5} 2^{12} \frac{1}{\{[256(51)]+1\}}$ (11)! $\left[\frac{1}{7!} - \frac{1}{5!} \frac{1}{(2n-1)^2 \pi^2} + \frac{1}{3!} \frac{1}{(2n-1)^4 \pi^4} - \frac{1}{(2n-1)^6 \pi^6}\right] =$

$\frac{4096(39916800)}{13057} \left[\frac{1}{5040} \frac{5\pi^5}{1536} - \frac{1}{120\pi^2} \frac{61\pi^7}{184320} + \frac{1}{6\pi^4} \frac{1385\pi^9}{41287680} - \frac{1}{\pi^6} \frac{50521\pi^{11}}{14863564800}\right] = \pi^5$

k = 4 gives series as $\sum_{n=1}^{\infty} \frac{(-1)^{n+1}}{(2n-1)^5} 2^{14} \frac{1}{\{[1024(74)]+1\}}$ (13)! $\left[\frac{1}{9!} - \frac{1}{7!} \frac{1}{(2n-1)^2 \pi^2} + \frac{1}{5!} \frac{1}{(2n-1)^4 \pi^4} - \frac{1}{3!} \frac{1}{(2n-1)^6 \pi^6} + \frac{1}{(2n-1)^8 \pi^8}\right] = \frac{16384(6227020800)}{75777} \left[\frac{1}{362880} \frac{5\pi^5}{1536} - \frac{1}{5040\pi^2} \frac{61\pi^7}{184320} + \frac{1}{120\pi^4} \frac{1385\pi^9}{41287680} - \frac{1}{6\pi^6} \frac{50521\pi^{11}}{14863564800} + \frac{1}{\pi^8} \frac{2702765\pi^{13}}{16384(479001600)}\right] = \pi^5$

On similar lines, the series representation (18) can be proved for any integer k value.

### 3. <u>Series representations for odd powers $\pi^2$, $\pi^4$ and $\pi^6$</u>

### 3a. The series for $\pi^2$

The series for $\pi^2$ is presented as

$$\pi^2 = \sum_{n=1}^{\infty} \frac{1}{n^2} 2(2k+3) \quad (2k+1)! \sum_{j=0}^{k} \left[-\frac{1}{(n)^2 \pi^2}\right]^j \frac{1}{(2k-2j+1)!} \qquad (19)$$

where j, k are integers and $\geq 0$.

### <u>Verification of the series for some k values</u>

k = 0 gives series as $\sum_{n=1}^{\infty}\frac{1}{n^2} 6 \ (1)! \ \left[-\frac{1}{(n)^2\pi^2}\right]^0 \frac{1}{(1)!} = 6 \cdot \frac{\pi^2}{6} = \pi^2$

k = 1 gives series as $\sum_{n=1}^{\infty}\frac{1}{n^2} 10 \ (3)! \ \left[\frac{1}{3!} - \frac{1}{(n)^2\pi^2}\right] = 60 \left[\frac{1}{6}\frac{\pi^2}{6} - \frac{1}{\pi^2}\frac{\pi^4}{90}\right] = \pi^2$

k = 2 gives series as $\sum_{n=1}^{\infty}\frac{1}{n^2} 14 \ (5)! \ \left[\frac{1}{5!} - \frac{1}{3!}\frac{1}{(n)^2\pi^2} + \frac{1}{(n)^4\pi^4}\right] = 1680 \left[\frac{1}{120}\frac{\pi^2}{6} - \frac{1}{6}\frac{1}{\pi^2}\frac{\pi^4}{90} + \frac{1}{\pi^4}\frac{\pi^6}{945}\right] = \pi^2$

k = 3 gives series as $\sum_{n=1}^{\infty}\frac{1}{n^2} 18 \ (7)! \ \left[\frac{1}{7!} - \frac{1}{5!}\frac{1}{(n)^2\pi^2} + \frac{1}{3!}\frac{1}{(n)^4\pi^4} - \frac{1}{(n)^6\pi^6}\right] =$

$90720 \left[\frac{1}{5040}\frac{\pi^2}{6} - \frac{1}{120}\frac{1}{\pi^2}\frac{\pi^4}{90} + \frac{1}{6}\frac{1}{\pi^4}\frac{\pi^6}{945} - \frac{1}{\pi^6}\frac{\pi^6}{9450}\right] = \pi^2$

k = 4 gives series as $\sum_{n=1}^{\infty}\frac{1}{n^2} 22 \ (9)! \ \left[\frac{1}{9!} - \frac{1}{7!}\frac{1}{(n)^2\pi^2} + \frac{1}{5!}\frac{1}{(n)^4\pi^4} - \frac{1}{3!}\frac{1}{(n)^6\pi^6} + \frac{1}{(n)^8\pi^8}\right] =$

$7983360\left[\frac{1}{362880}\frac{\pi^2}{6} - \frac{1}{5040}\frac{1}{\pi^2}\frac{\pi^4}{90} + \frac{1}{120}\frac{1}{\pi^4}\frac{\pi^6}{945} - \frac{1}{6}\frac{1}{\pi^6}\frac{\pi^8}{9450} + \frac{1}{\pi^6}\frac{\pi^8}{93555}\right] = \pi^2$

On similar lines, the series representation (19) can be proved for any integer k value.

### 3b. The series for $\pi^4$

The series for $\pi^4$ is presented as

$$\pi^4 = \sum_{n=1}^{\infty}\frac{1}{n^4} 6(2k+5)(2k+3) \ (2k+1)! \ \sum_{j=0}^{k}\left[-\frac{1}{(n)^2\pi^2}\right]^j \frac{1}{(2k-2j+1)!} \qquad (20)$$

where j, k are integers and $\geq 0$.

### Verification of the series for some k values

k = 0 gives series as $\sum_{n=1}^{\infty}\frac{1}{n^4} 6\ (5)(3) \ (1)! \ \left[-\frac{1}{(n)^2\pi^2}\right]^0 \frac{1}{(1)!} = 90 \cdot \frac{\pi^4}{90} = \pi^4$

k = 1 gives series as $\sum_{n=1}^{\infty}\frac{1}{n^4} 6\ (7)(5) \ (3)! \ \left[\frac{1}{3!} - \frac{1}{(n)^2\pi^2}\right] = 1260 \left[\frac{1}{6}\frac{\pi^4}{90} - \frac{1}{\pi^2}\frac{\pi^6}{945}\right] = \pi^4$

k = 2 gives series as $\sum_{n=1}^{\infty} \frac{1}{n^4} 6(9)(7) \ (5)! \left[ \frac{1}{5!} - \frac{1}{3!} \frac{1}{(n)^2 \pi^2} + \frac{1}{(n)^4 \pi^4} \right] =$

$45360 \left[ \frac{1}{120} \frac{\pi^4}{90} - \frac{1}{6} \frac{1}{\pi^2} \frac{\pi^6}{945} + \frac{1}{\pi^4} \frac{\pi^8}{9450} \right] = \pi^4$

k = 3 gives series as $\sum_{n=1}^{\infty} \frac{1}{n^4} 6(11)(9) \ (7)! \left[ \frac{1}{7!} - \frac{1}{5!} \frac{1}{(n)^2 \pi^2} + \frac{1}{3!} \frac{1}{(n)^4 \pi^4} - \frac{1}{(n)^6 \pi^6} \right] =$

$2993760 \left[ \frac{1}{5040} \frac{\pi^4}{90} - \frac{1}{120} \frac{1}{\pi^2} \frac{\pi^6}{945} + \frac{1}{6} \frac{1}{\pi^4} \frac{\pi^8}{9450} - \frac{1}{\pi^6} \frac{\pi^{10}}{93555} \right] = \pi^4$

k = 4 gives series as $\sum_{n=1}^{\infty} \frac{1}{n^4} 6(13)(11) \ (9)! \left[ \frac{1}{9!} - \frac{1}{7!} \frac{1}{(n)^2 \pi^2} + \frac{1}{5!} \frac{1}{(n)^4 \pi^4} - \frac{1}{3!} \frac{1}{(n)^6 \pi^6} + \frac{1}{(n)^8 \pi^8} \right] =$

$311351040 \left[ \frac{1}{362880} \frac{\pi^4}{90} - \frac{1}{5040} \frac{1}{\pi^2} \frac{\pi^6}{945} + \frac{1}{120} \frac{1}{\pi^4} \frac{\pi^8}{9450} - \frac{1}{6} \frac{1}{\pi^6} \frac{\pi^{10}}{93555} + \frac{1}{\pi^8} \frac{691\pi^{12}}{638512875} \right] = \pi^4$

On similar lines, the series representation (20) can be proved for any integer k value.

### 3c. The series for $\pi^6$

The series for $\pi^6$ is presented as

$$\pi^6 = \sum_{n=1}^{\infty} \frac{1}{n^6} \frac{45(2k+7)(2k+5)(2k+3)}{(k+5)} \ (2k+1)! \ \sum_{j=0}^{k} \left[ -\frac{1}{(n)^2 \pi^2} \right]^j \frac{1}{(2k-2j+1)!} \qquad (21)$$

where j, k are integers and $\geq 0$.

### Verification of the series for some k values

k = 0 gives series as $\sum_{n=1}^{\infty} \frac{1}{n^6} 45(7)(3) \ (1)! \left[ -\frac{1}{(n)^2 \pi^2} \right]^0 \frac{1}{(1)!} = 945 \cdot \frac{\pi^6}{945} = \pi^6$

k = 1 gives series as $\sum_{n=1}^{\infty} \frac{1}{n^6} 45(9)(7)(5)\frac{1}{6} \ (3)! \left[ \frac{1}{3!} - \frac{1}{(n)^2 \pi^2} \right] = 14175 \left[ \frac{1}{6} \frac{\pi^6}{945} - \frac{1}{\pi^2} \frac{\pi^8}{9450} \right] = \pi^6$

k = 2 gives series as $\sum_{n=1}^{\infty} \frac{1}{n^6} 45(11)(9) \ (5)! \left[ \frac{1}{5!} - \frac{1}{3!} \frac{1}{(n)^2 \pi^2} + \frac{1}{(n)^4 \pi^4} \right] =$

$534600 \left[ \frac{1}{120} \frac{\pi^6}{945} - \frac{1}{6} \frac{1}{\pi^2} \frac{\pi^8}{9450} + \frac{1}{\pi^4} \frac{\pi^{10}}{93555} \right] = \pi^6$

k = 3 gives series as $\sum_{n=1}^{\infty} \frac{1}{n^6}$ 45 (13)(11)(9)$\frac{1}{8}$ (7)! $\left[ \frac{1}{7!} - \frac{1}{5!} \frac{1}{(n)^2 \pi^2} + \frac{1}{3!} \frac{1}{(n)^4 \pi^4} - \frac{1}{(n)^6 \pi^6} \right]$ = $36486450 \left[ \frac{1}{5040} \frac{\pi^6}{945} - \frac{1}{120} \frac{1}{\pi^2} \frac{\pi^8}{9450} + \frac{1}{6} \frac{1}{\pi^4} \frac{\pi^{10}}{93555} - \frac{1}{\pi^6} \frac{691\pi^{12}}{638512875} \right] = \pi^6$

k = 4 gives series as $\sum_{n=1}^{\infty} \frac{1}{n^6}$ 45(15)(13)(11)$\frac{1}{9}$ (9)! $\left[ \frac{1}{9!} - \frac{1}{7!} \frac{1}{(n)^2 \pi^2} + \frac{1}{5!} \frac{1}{(n)^4 \pi^4} - \frac{1}{3!} \frac{1}{(n)^6 \pi^6} + \frac{1}{(n)^8 \pi^8} \right]$ = $3891888000 \left[ \frac{1}{362880} \frac{\pi^6}{945} - \frac{1}{5040} \frac{1}{\pi^2} \frac{\pi^8}{9450} + \frac{1}{120} \frac{1}{\pi^4} \frac{\pi^{10}}{93555} - \frac{1}{6} \frac{1}{\pi^6} \frac{691\pi^{12}}{638512875} + \frac{1}{\pi^8} \frac{2\pi^{14}}{18243225} \right] = \pi^6$

On similar lines, the series representation (21) can be proved for any integer k value.